\input ./style/arxiv-general.cfg
\documentclass[MSNbibl,number,citesort,seceqn,dvips]{arxbj}
\makeatletter
   \@ifpackageloaded{graphicx}{}{\usepackage{graphicx}}
\makeatother
\usepackage{upgreek}
\usepackage[small,nohug,heads=littlevee]{diagrams}

\aid{0}
\volume{21}
\issue{3}
\pubyear{2015}
\firstpage{1289}
\lastpage{1303}
\doi{10.3150/13-BEJ541} 

\makeatletter
\newcommand{\rrvert}{\vert}
\newcommand{\llvert}{\vert}

\def\shalf{\mbox{{\tiny$\frac{1}{2}$}}}
\def\P{\mathrm{P}}

\newtheorem{theorem}{Theorem}[section]

\def\halft{\frac{1}{2t}}
\newcommand{\fF}{\ensuremath{\mathfrak{F}}}
\newcommand{\cA}{\mathcal{A}}
\newcommand{\cM}{\mathcal{M}}
\newcommand{\cR}{\mathcal{R}}
\newcommand{\dd}{\mathrm{d}}
\def\E{\mathrm{E}}
\def\var{\operatorname{var}}
\def\cov{\operatorname{cov}}
\makeatother

\begin{document}
\begin{frontmatter}

\title{Applying Dynkin's isomorphism: An~alternative approach to
understand the~Markov property of the de Wijs process}
\runtitle{Markov property of the de Wijs process}

\begin{aug}
\author{\inits{D.}\fnms{Debashis} \snm{Mondal}\corref{}\ead[label=e1]{debashis@stat.oregonstate.edu}}
\address{Department of Statistics, Oregon State University, 44 Kidder Hall,
Corvallis, OR 97331, USA.\\ \printead{e1}}
\end{aug}

\received{\smonth{1} \syear{2010}}
\revised{\smonth{5} \syear{2013}}

%
\begin{abstract}
Dynkin's (\textit{Bull. Amer. Math. Soc.} \textbf{3} (1980) 975--999)
seminal work associates a multidimensional transient
symmetric Markov
process with a multidimensional Gaussian random field. This
association, known as Dynkin's isomorphism, has profoundly
influenced the studies of Markov properties of generalized Gaussian
random fields. Extending Dykin's isomorphism, we study here a
particular generalized Gaussian Markov random field, namely, the de
Wijs process that originated in Georges Matheron's pioneering work on
mining geostatistics and, following McCullagh (\textit{Ann. Statist.}
\textbf{30} (2002) 1225--1310), is now
receiving renewed attention in spatial statistics. This extension of
Dynkin's theory associates the de Wijs process with the (recurrent)
Brownian motion on the two dimensional plane, grants us further insight
into Matheron's kriging formula
for the de Wijs process and highlights previously unexplored
relationships of the central Markov models in spatial statistics with
Markov processes on the plane.
\end{abstract}

%

\end{frontmatter}

\section{Introduction} \label{sec:introduction}

After originating in the pioneering work of Georges Matheron, the de Wijs
process enjoyed a significant and extensive role in early
geostatistical literature \cite{r6,r26,r27a,r18}. McCullagh's
\cite{r29} recent work has revived interest in the de Wijs process,
both theoretically and in a growing range of
applications in spatial statistics; see, for example, \cite
{r4,r33,r5,r8,r9,r10,r30,r34,r35,n1}. In particular, Mondal \cite{r33} and
 Besag and Mondal
\cite{r5}
established a connection between Gaussian Markov random fields on
two-dimensional lattices
and the de Wijs process on the Euclidean plane, which emerges as a
scaling limit of the former.
 McCullagh and Clifford
 \cite{r30} analyzed
agricultural uniformity trials using a spatial formulation that is
based on the de Wijs process and a Gaussian white noise random
field. See also the related work by Clifford
\cite{r8,r9}, Clifford et al. \cite{r10}. Mondal \cite{r34}
considers the
exponential functional of the de Wijs
process to construct a generalized Cox process to study disease
mappings. Mondal \cite{r35} indicates a link between the de Wijs process
and Tobler's \cite{r51} pycnoplylectic interpolation based on the Laplace
splines. Dutta and Mondal \cite{n1} make explicit
use of the connection
between intrinsic autoregressions and the de Wijs process and provide
approximate matrix free computations for residual maximum likelihood
methods for the latter. Furthermore, outside the statistics literature,
the de Wijs process appears to originate separately in quantum physics
and statistical mechanics as the massless case of the free Gaussian
field; see, for example, Chapters~6 and 7 of \cite{r16},
and in recent probability literature, this massless case has become a
subject of intense study; see, for example, \cite{r20,r47}.

Technically, the de Wijs process is a generalized Gaussian random
field (\cite{r15}, Chapter III), whose index set
is a certain class of contrasts, that is, non-atomic signed Borel measures
on the Euclidean plane with zero total mass. This process
corresponds to the logarithmic variogram model and is a
generalization of the Brownian motion in two dimensions. It acquires
Markov and conformal invariance properties \cite{r29,r30} and is
first-order intrinsic in the
sense of \cite{r54} and \cite{r12}. The Markov property of
the de Wijs process was already known to Matheron \cite{r27,r27a}, who
viewed it from the perspective of kriging predictions. Consider a
mean zero Gaussian random field $\{ U(x)\dvtx  x \in\cR^2\}$.
Let $B$ denote a closed contour in $\cR^2$. It is natural to
call the random field \emph{Markov} if its values along the curve
$B$ determine the ordinary kriging predictor for the value $U(x_0)$
at a point $x_0$ in the interior of $B$, given its values on and in
the exterior of $B$. This Markovian characterization leads to a
kriging predictor for the random variable $U(x_0)$ that takes the
form of a contour integral
\[
\E \bigl( U(x_0) \mid U(x), x \in B \bigr) = \int
_B v(x, x_0) U(x) \,\mathrm{d} x,
\]
where the coefficient function $v(x, x_0)$, $x \in B$ is such that
%
\begin{equation}
\label{eq:matheronkriging} \int_B v(x, x_0) \,\mathrm{d} x =1, \qquad \int
_B v(x, x_0) \cov\bigl(U(x), U
\bigl(x'\bigr)\bigr) \,\mathrm{d} x = \cov\bigl(U(x_0), U
\bigl(x'\bigr)\bigr)
\end{equation}
for every point $x'$ on $B$. When $\{ U(x) \dvtx  x \in\cR^2
\}$ is stationary and isotropic, Matheron \cite{r27}
deduces that it is
Markov if and only if
%
\begin{equation}
\label{eq:covgou} \cov\bigl(U(x), U(0)\bigr) \propto K_0 \bigl(a \|x\|\bigr),
\end{equation}
where $K_0$ is the Bessel function of order zero and $a$ is a positive constant.
Matheron further notes that his derivations remain intact for the
limiting case $a \downarrow0$ that corresponds to the logarithmic
covariance, $\cov(U(x), U(0)) = - \log( \|x \|)$ and thus provides
the Markovian characterization of the de Wijs process. However, the
respective random fields exist as generalized processes only, and the
details of a formal argument go beyond the above kriging formula.

This paper calls attention to the work of Dynkin \cite
{r13} to present a
mathematical formalism to
describe the above kriging formula of the de Wijs process, and to
connect the field of spatial statistics to the vast, and hitherto
unutilized, probabilistic literature on the Markov property of
generalized random fields. This body of literature constitutes a
fascinating part of probability, much of which emerged in the wake of
\cite{r24} and \cite{r31}. This corpus notably includes
\cite{r52,r53,r32,r37,r38,r19,r43,r44,r14,r21,r13,r17,r41,r42,r25,r46} and the references
therein. In many of these works, several notions of Markovianity
for generalized Gaussian random fields have emerged and their
interrelations and their connections to various related concepts often
form a good part
of their understanding. For example, a homogeneous
and isotropic generalized Gaussian random field whose spectral
density is inversely proportional to an even polynomial of the
frequencies satisfies a Markov property in
the sense of Holley and Stroock but may not be Markov in the sense of
Wong. It is interesting to note that Nelson's
\cite{r37} construction of the free Markov field on the plane actually
corresponds to the generalized Gaussian random field with covariance
given in (\ref{eq:covgou}). Using a slightly different notion of
Markovianity, Wong \cite{r53} arrives, much earlier
than Nelson, at the
conclusion that the only generalized Gaussian Markov random field
again has covariance (\ref{eq:covgou}). Kallianpur and Mandrekar \cite
{r19}, on the other hand, investigate Markov properties of a
generalized Gaussian random field in conjunction with its dual
random field. Ekhaguere \cite{r14} later provides links
between the Markov
property due to Nelson \cite{r37} and that due to
Wong \cite{r53}. In contrast,
Dynkin's \cite{r13} famous work marks an important departure
from these
earlier studies. In his study, covariances of a generalized Gaussian
random field are assumed to arise from the Green function of a
symmetric multidimensional Markov process, and the Markov property
of this generalized Gaussian random field is then derived from
the path properties of the multidimensional Markov process. Thus, for
example, the
Markov property of (\ref{eq:covgou}) can be understood in the context
of the
Markov property of an exponentially killed Brownian motion on the
plane. Here our focus is on the limiting case, namely, the de Wijs
process. Although its Markov property can be investigated using the
work of Nelson or Wong (e.g., by modifying Theorem~1.5 of \cite{r17} or by including the case $\alpha=0$ in Wong's
\cite{r53}
Theorem~2), we take up Dynkin's approach primarily
because it provides a precise and computable description of the
boundary condition in the kriging formula, and connects closely with
Matheron's work. We also piece together many scattered results and
extend some known ones to provide this new addition to the body of
literature that, respectively, followed Dynkin's and Matheron's works.

The remainder of the paper is structured as follows. Section~\ref{sec:dewijs}
introduces the de Wijs process as a homogeneous, isotropic and self-similar
generalized Gaussian random field. Section~\ref{sec:brownain} explores the association
of this process with Brownian motion. Here we show that the covariance
formula of the de Wijs process can be written explicitly in terms of an
additive function of the Brownian motion. Section~\ref{sec:markov} studies the Markov
property of the de Wijs process by extending the work of Dynkin \cite{r13}.
Here our main result, namely Theorem~\ref{th:markov}, provides a new
interpretation of Matheron's kriging formula in terms of the hitting
probabilities of the Brownian motion and as a generalization of the
Dirichlet problem. Section~\ref{sec:discussion} focuses on the practical relevance of
Matheron's kriging formula. It also considers the relevance of Dynkyn's
isomorphisms in lattice approximations of the de Wijs process and
concludes with a discussion on the screening effect in kriging.

\section{De Wijs process} \label{sec:dewijs}

In this paper, a \emph{generalized random field} on the Euclidean
plane $\mathcal{R}^2$ is a stochastic process $\{ Z_\sigma \dvtx \sigma
\in\cM\}$ indexed by a vector space $\cM$ of non-atomic signed
Borel measures on the plane that have total mass zero. We view $Z$ to
be a linear functional from the vector space $\cM$ to the real numbers
such that
\[
Z_{b \sigma+ \mathrm{d} \nu} = b Z_\sigma+ \mathrm{d} Z_\nu\qquad  \mbox{for all }
\sigma , \nu\in\cM, \mbox{ and for all } b, d \in\cR.
\]
We think of the random variable $Z_\sigma$ as a spatial \emph
{contrast}; for instance, if two plots have unit area and $\sigma$
has a Lebesgue density that is proportional to the difference of the
respective indicator functions, then $Z_\sigma$ might represent the
difference of crop yields on these plots.

The generalized random field $\{ Z_\sigma\dvtx  \sigma\in\cM\}$ is
said to be \emph{homogeneous} if its distribution remains invariant to
planar translations, and \emph{isotropic} if its distribution remains
invariant to planar rotations. Furthermore, such a homogeneous
isotropic generalized random field is \emph{Gaussian}
if all finite dimensional marginal distributions are multivariate
normal with
$\E Z_\sigma= 0$ and
%
\begin{equation}
\label{eq:cov} \cov( Z_\sigma, Z_\nu) = - \int\!\!\!\int\varphi\bigl(
\|x-y\|\bigr) \sigma(\dd x) \nu(\dd y)
\end{equation}
for some real-valued function $\varphi$ and all non-atomic signed
measures $\sigma, \nu\in\cM$. Note that in the above covariance
formula we can add to $\varphi(\|x-y\|)$ a function $f_1(x) + f_2(y)$
without affecting the integral, and so $\varphi$ actually belongs to a
suitable quotient space (modulo the infinite dimensional subspace of
additive functions). In subsequent discussions, we will be implicit
about this equivalence relation in the description of $\varphi$.
The function
\[
c(x) = \varphi\bigl(\| x \|\bigr),\qquad  x \in\mathcal{R}^2
\]
is then called the \emph{generalized variogram} or $-c(x)$ the \emph
{generalized covariance
function} of the generalized random field. Let $\hat\sigma$ and
$\hat\nu$ denote the Fourier transforms of $\sigma$ and $\nu$. We
can then write (\ref{eq:cov}) as
%
\begin{equation}
\label{eq:spectral} \cov( Z_\sigma, Z_\nu) = \int\hat{\sigma}(x)
\overline{\hat{\nu}(x)} S(\dd x)
\end{equation}
for a certain non-negative tempered measure $S$ which is called the
\emph{spectral measure} of the generalized random field (\cite{r15}, page 264). If the spectral measure is absolutely
continuous with respect to the Lebesgue measure on $\mathcal{R}^2$,
its Lebesgue density $s(x)$, $x \in\mathcal{R}^2$ is called the
\emph{spectral density}. Under slight regularity conditions, the
generalized covariance function $c$ and the spectral density $s$ are
Fourier transforms of each other.

Specifically, consider $\cM$ to be the space of signed Borel
measures $\sigma$ on the Euclidean plane $\mathcal{R}^2$ that
satisfy
%
\begin{equation}
\label{eq:finite} \int\!\!\!\int\bigl\llvert \log \bigl( \| x-y \| \bigr) \bigr\rrvert |\sigma|(\dd
x) |\sigma| (\dd y) < \infty
\end{equation}
and have total mass zero. The \emph{de Wijs process} is then the
homogeneous, isotropic and self-similar generalized Gaussian random
field (i.e., its distribution also remains invariant to changes of
scale) on
$\mathcal{R}^2$ with index set $\cM$ such that $\E Z_\sigma= 0$
and
%
\begin{equation}
\label{eq:deWijs} \cov( Z_\sigma, Z_\nu) = \langle\sigma, \nu
\rangle_{\cM} = - \int\!\!\!\int \log \bigl( \| x-y \| \bigr) \sigma(\dd x) \nu(\dd y)
\end{equation}
for all signed measures $\sigma, \nu\in\cM$. Note that $\cM$
has an inner product space structure with inner product $\langle\sigma,
\nu\rangle_{\cM}$ and norm
\[
\| \sigma\|_{\cM} = \langle\sigma, \sigma\rangle_{\cM}^{1/2}.
\]
Indeed, Corollary~2.5 of \cite{r28} implies
that $\cM$ is a
vector space, and by Corollary~2.4 and Remark~3.3 in the same
reference $\| \sigma\|_{\cM} \geq0$ with equality if and only if
$\sigma$ vanishes identically. The positive definiteness of the
covariance matrices associated with the de Wijs structure
(\ref{eq:deWijs}) is an immediate consequence of the Gram matrix
property. Thus, the de Wijs process has logarithmic variogram; that is, the
representation (\ref{eq:cov}) holds with the generalized variogram
function $c(x) = \varphi(\|x\|) = \log\|x\|$ for $x \in
\mathcal{R}^2$, and its spectral density is
\[
s(x) = \frac{1}{2 \uppi\| x \|^2}, \qquad x \in\mathcal{R}^2.
\]

\section{Association with Brownian motion} \label{sec:brownain}

We set $T=[0,\infty)$ for consistency in what follows and let $\{ W_t,
t \in T\}$ be the Brownian motion on the two-dimensional Euclidean
plane. Thus, with probability 1, the function $t \rightarrow W(t)$ is
continuous in $t$, the components of the increment $W_{t+u} -W_u$ are
independent Gaussian random variables each with mean $0$ and variance
$t$, and the process $\{ W_t, t \in T\}$ has stationary and independent
increments. For every $x$ on the plane, let $\P_x$ denote the
probability law of
$\{ W_t, t \in T \}$ starting at $x$ and let $\E_x$ be its expectation
under $\P_x$. For every $t \in T$, let the
sub-$\sigma$-field $\fF_t$ consist of events observable up to time
$t$, which is the minimum $\sigma$-field generated by $\{W_u\dvtx
0 \le u \le t\}$. Define $\fF_\infty$ to be the minimum $\sigma
$-field containing $\bigcup_{t \in T} \fF_t$. The Markov property of the
Brownian motion asserts that the conditional law of $\{W_t, t\ge u\}$
given $\{W_t, 0\le t \le u\}$ depends on $W$ only through $W_u$. In
other words,
\[
\E_x (F J) = \E_x (F \mathrm{E}_{W_u} J)
\]
for every $x$ on the plane, every $\fF_u$ measurable positive function
$F$ and every measurable function $J$ that depend only on $\{W_t, t\ge
u\}$. In particular, the expectation is calculated first with respect
to the conditional law of $\{W_t, t\ge u\}$ given $\{W_t, 0\le t \le u\}
$, and then with respect to the marginal law of $\{W_t, 0\le t \le u\}
$. An important generalization of the Markov property is the strong
Markov property. When $\tau$ is a stopping time, define the stopping
field $\fF_\tau$ to be the $\sigma$-algebra consisting of all events
$A \in\fF_\infty$ such that $A \cap\{ \tau\le t\} \in\fF_t$ for
every $t \ge0$. Then the strong Markov property implies that
\[
\E_x (F J) = \E_x (F \mathrm{E}_{W_\tau} J)
\]
for every $x$ on the plane, every stopping time $\tau$, every $\fF
_\tau$ measurable positive function $F$ and every measurable function
$J$ that depends on $\{W_t, t \ge\tau\}$.

Next we define the Green function of the Brownian motion $\{ W_t, t \in
T\}$. Typically, the Green function is defined for a transient Markov
process as the time integral of its transition probability density
function, and Dynkin's theory is essentially based on the fact that the
Green function of a transient symmetric Markov process can be
interpreted as the covariance of a centered Gaussian process. However,
the Brownian motion on the plane is recurrent and hence its transition
probability density function
\[
p_t(x,y) = (2 \uppi t)^{-1} \exp\biggl\{ - \halft\| y-x
\|^2\biggr\}
\]
is not integrable with respect to $t$. Thus, we need a modification
that will allow us to define the Green function of the Brownian motion
and extend Dynkin's result in a straightforward fashion. To this end,
we fix a point $x_0$ on the unit circle and consider $q_t(x,y) =
p_t(0,y-x) - p_t(0,x_0)$. We then apply the definition of Port and Stone
(\cite{r39}, page 70) and obtain the \emph{Green function} or
\emph{the
potential kernel} of $\{ W_t, t \in T \}$ as
%
\begin{equation}
\label{eq:disomorphism} g(x,y) = \int_0^\infty
q_t(x,y) \,\mathrm{d}t = -\uppi^{-1} \log\| y-x \|.
\end{equation}
The choice of $x_0$ will not matter, as we shall see in what follows.

However, note that the covariances of the de Wijs process now satisfy
the relationship
\[
\cov(Z_\sigma, Z_\nu) = \langle\sigma, \nu
\rangle_{\cM} = \int\!\!\!\int\uppi g(x,y) \sigma(\mathrm{d} x) \nu(\mathrm{d} y).
\]
Thus, we say that the above relationship \emph{associates} the de Wijs
process with
the Brownian motion $\{W_t, t\in T\}$, opening an avenue
for exploring the properties of the former from those of the latter.
For every $\nu\in\cM$, we now get
\[
\int q_t(x,y) \nu(\mathrm{d}y) = \int p_t(x,y) \nu(\mathrm{d}y) -
p_t\bigl(0, \|x_0\|\bigr) \int\nu(\mathrm{d}y) = \int p_t(x,y)
\nu(\mathrm{d}y)
\]
and therefore
\[
\int g(x,y) \nu(\mathrm{d}y) = \int\!\!\!\int p_t(x,y) \nu(\mathrm{d}y)\, \mathrm{d}t.
\]
In other words, the term involving $x_0$ disappears from the right-hand
side of the previous equation. Consequently, when $\nu$ is absolutely
continuous with the Radon--Nykodyn derivative
$\nu(\mathrm{d} y) = \rho(y) \,\mathrm{d}y$, the above equation becomes
\[
\int g(x,y) \nu(\mathrm{d}y) = \E_x \int\rho(W_t) \,\mathrm{d}t.
\]
Now equation (\ref{eq:disomorphism}) can be identified with
$\langle\sigma, \nu\rangle_{\cM} = \uppi\E_\sigma\int\rho(W_t) \,\mathrm{d}t$,
where $\E_\sigma$ is the expectation under the probability law of $\{
W_t, t \in T \}$ with initial signed measure $\sigma$; that is,
\[
\E_\sigma\int\rho(W_t) \,\mathrm{d}t = \int\E_x \int
\rho(W_t) \,\mathrm{d}t \,\sigma(\mathrm{d}x).
\]
We can thus define the \emph{additive function} of the Brownian motion
by the measure
%
\begin{equation}
{\label{eq:additivef}} A_\nu(Q) = \int_Q
\rho(W_t ) \,\mathrm{d}t
\end{equation}
that satisfies the property that, for every interval $I = (s,u)$ with
$s<u$, $A_\nu(I)$ is
a functional of $\{ W_t, t \in I\}$, and
%
\begin{equation}
{\label{eq:addfandcov}} \cov(Z_\sigma, Z_\nu) = \langle\sigma, \nu
\rangle_{\cM} = \uppi \E_\sigma A_\nu(T).
\end{equation}
The collection of all signed measures $\nu\in\cM$ that are
absolutely continuous forms a dense subspace of $\cM$. By passage to
limit, it then follows that for every $\nu\in\cM$ there
exists an additive functional of the Brownian motion such that
the above equation holds. In addition, the strong Markov property of
the Brownian motion
takes the following form
%
\begin{equation}
\label{eq:strongmarkov} \E_\sigma F A_\nu(\tau+ Q) = \E_\sigma
F \mathrm{E}_{W_\tau} A_\nu(Q)
\end{equation}
for every $\sigma, \nu\in\cM$, for every Borel subset $Q$ of $T$,
and for every $\tau$ and $F$ as defined earlier. The strengthened
relationship that emerges from equation (\ref{eq:addfandcov}) in
conjunction with equation (\ref{eq:strongmarkov}) now paves the way to
use the Brownian paths to study the properties of the de Wijs process.

\section{Markov property of the de Wijs process}\label{sec:markov}

As in Section~\ref{sec:introduction}, take $B$ to be a simple closed contour on the plane.
Then $B$ divides the entire plane into two components, namely, the
bounded interior and the unbounded exterior. Let $B_{\mathrm{I}}$ denote the open
interior of $B$ with closure $\bar B_{\mathrm{I}}$. Similarly, let $B_\mathrm{E}$ be the
open exterior of $B$ with closure $\bar B_{\mathrm{E}}$. Our first task is to
describe the values of the de Wijs process on the boundary $B$, and on
the inside and the outside of $B$ (e.g., on sets $B_\mathrm{I}$ and $\bar B_{\mathrm{E}}$).
To this end, there are two approaches. The first approach is due to
 \cite{r38}. Here we describe the values of the
de Wijs process on
an open set $G$ by the minimum sigma field $\cA_G$ generated by all
$Z_\sigma$ such that $\sigma\in\cM$ and support of $\sigma$ is
compactly contained in $G$. Then, for any closed set $C$ the values of
the de Wijs process is described by the sigma field
\[
\cA_C = \bigcap_{G \supset C}
\cA_{G},
\]
where the intersection is taken over all open sets $G$ that contain
$C$. Thus, $\cA_{B_\mathrm{I}}$, $\cA_{\bar B_{\mathrm{E}}}$ and $\cA_B$ represent the
values of the de Wijs process on the inside, outside and on the
boundary, respectively, and the Markov property of the de Wijs process
asserts that, for any $\sigma\in\cM$ with support of $\sigma$
compactly contained in $B_\mathrm{I}$,
\[
\E ( Z_\sigma| \cA_{\bar B_{\mathrm{E}}} ) = \E ( Z_\sigma|
\cA_{B} )
\]
almost surely in the probability distribution of the de Wijs process.
Note that, by construction, the minimum sigma fields $\cA_{\bar B_{\mathrm{E}}}$
and $\cA_B$ contains neighborhood information, not just the
information on the set. The second approach adopted by Dynkin
\cite{r13} is
a simplified version of the above and goes as follows. For a close set
$C$, we define $\cM_C$ to be the set of signed Borel measures $\sigma
\in\cM$ that do not charge on its complement. Then, following
\cite{r13}, the minimum sigma field generated
by the collection of random
variables $\{ Z_\sigma\dvtx  \sigma\in\cM_C\}$ describes the values of
the de Wijs process on the $C$. The Markov property is described in the
usual fashion, namely, the values on the inside ($\bar B_{\mathrm{I}}$) and the
outside $(\bar B_{\mathrm{E}})$ of $B$ are conditionally independent given the
values on $B$. This leads to the following theorem.
%
\begin{theorem}\label{th:markov}
Let $D$ be any closed set on the plane with a simply connected open
interior, and let $\tau_D$ be the first hitting time of $D$ by the
Brownian motion $\{W_t, t\in T\}$ starting at $x$ on the plane. Denote
by $V_x$ the probability measure of $W_{\tau_D}$ conditioned on
$W_0=x$. For any initial signed measure $\sigma\in\cM$ which does
not charge on the set $D$, identify $\sigma_D$ with the signed measure
induced by $W_t$ at the first hitting time of $D$; that is,
\[
\sigma_D(G) = \int\P_x( W_{\tau_D} \in G)
\sigma(\mathrm{d}x) = \int V_x(G) \sigma(\mathrm{d}x)
\]
for all Borel subsets $G$ on the plane. Then the conditional
expectation of $Z_\sigma$, given the values of the de Wijs process on
$D$, is identical to $Z_{\sigma_D}$. In other words,
%
\begin{equation}
\label{eq:dynkinkriging} \E \bigl(Z_\sigma| \{Z_\nu, \nu\in
\cM_D\} \bigr) = Z_{\sigma_D}.
\end{equation}
\end{theorem}

Before we turn to the proof, let us first see how Matheron's kriging
formula enters into the above theorem. Equation (\ref
{eq:dynkinkriging}) implies that for all $\sigma\in\cM$ and $\nu\in
\cM_D$
\[
\E Z_\sigma Z_\nu= \E Z_{\sigma_D} Z_\nu.
\]
Hence, formula (\ref{eq:deWijs}) applies. This along with the
definition of $\sigma_D$ produces the identity
\[
\int\!\!\!\int\log \bigl( \| x-y \| \bigr) \sigma(\dd x) \nu(\dd y) = \int\!\!\!\int\int\log \bigl( \|
x'-y \| \bigr) V_{x} \bigl(\dd x'\bigr)
\sigma(\dd x) \nu(\dd y)
\]
for all $\sigma\in\cM$ and $\nu\in\cM_D$. Consequently, if $x$ is
in the interior of boundary $B$, we have
\[
\int_B V_{x} \bigl(\dd x'\bigr)
=1, \qquad \log \bigl( \| x-y \| \bigr) = \int_{B} \log \bigl( \|
x'-y \| \bigr) V_{x} \bigl(\dd x'\bigr)
\]
for every point $y$ on $B$. In short, the coefficient function
$v(x',x_0)$ in Matheron's kriging formula (\ref{eq:matheronkriging})
is the derivative of $V_{x_0}$ at $x'$, and thus corresponds to the
probability density function of the Brownian motion at the first
hitting time $\tau_D$, a crucial fact that has arguably been missing
from the geostatistitical literature. Furthermore, in order for $Z$ to
be a linear functional from the vector space $\cM$ to the real
numbers, we can imagine $Z_\sigma$ as an integral of the form $\int
Z_x \sigma(\dd x)$, where the notation $Z_x$ suggests a point-wise
intrinsic process with $\var(Z_x-Z_{x'}) = -\log ( \| x-x' \|
)$. This very imagination of a point-wise $Z_x$ allows us to
describe $\hat{Z}_{x}= \int_B Z_{x'} V_{x} (\dd x')$ as the kriged
value of $Z_{x}$, for an $x$ in the interior of $B$. We can take this
point further, and even describe the kriging formula from a different
angle. First, let $\nabla$ denote the Laplace operator on the plane.
If $\nu$ is twice differentiable, Theorem~3 of \cite{r40},
page 525,
implies
%
\begin{equation}
\label{eq:laplace} -2 \uppi\langle\sigma, \nabla\nu\rangle_{\cM} = \int
\sigma(x) \nu(x) \,\dd x,
\end{equation}
which surprisingly asserts that $\nabla Z_x$ and $Z_{x'}$ behave as two
mean zero uncorrelated Gaussian random variables for $x \neq x'$, and
in turn, suggests that the kriged values of $\nabla Z_x$ on the
interior of $B$ are all zero. An interchange of the Laplace operator
and the conditional expectation on $Z_x$ then produce
\[
\nabla\hat{Z}_x =0
\]
on the interior of $B$, implying that the kriging problem is a
generalization of the Dirichlet problem in mathematics. Indeed, the
probability literature reaffirms that the boundary values of the
Brownian motion at first hitting time solve the standard version of the
Dirichlet problem and therefore the Matheron's kriging formula (\ref
{eq:matheronkriging}) can be seen as a generalization. For significance
of the Dirichlet problem in recent spatial statistics, we refer to the
discussion in \cite{r29}. We now return to the
proof of the theorem.
\begin{pf*}{Proof of Theorem~\ref{th:markov}}
First, we verify that $\sigma_D \in\cM_D$. As a first step, we argue
that $\sigma_D$ belongs to $\cM$. Since $\int\sigma_D(\mathrm{d}x) = \int
\sigma(\dd x) = 0$, $\sigma_D$ represents a signed Borel measure with
total mass zero. Now, for a non-negative measure $\mu$ for which the
integral $h_\mu(x) = \int g(x,y) \mu(\dd y)$ is finite for every $x$,
the results of \cite{r7}, pages 193--194, give the identity
\[
h_\mu(x) - \E_x h_\mu(W_t) = \int
g(x,y) \mu(\mathrm{d}y) - \E_x \int g(W_t,y) \mu(\mathrm{d}y) = \int\!\!\!\int
_0^t p_s(x,y) \,\mathrm{d}s \,\mu(\mathrm{d}y) \ge0.
\]
Consequently, $h_\mu$ defines an \textit{excessive measure}, and
$h_\mu(x) \ge \E_x h_\mu(W_t)$. Hence, the choices $\mu= \nu^+$ and
$\mu=\nu^-$ yield
\[
h_{\nu^+}(x) \ge\E_x h_{\nu^+} (W_t),\qquad
h_{\nu^-}(x) \ge\E_x h_{\nu^-} (W_t).
\]
It then follows that
\[
\bigl\langle\sigma_D^+, \nu^+\bigr\rangle_{\cM} = \uppi
\E_{\sigma^+} h_{\nu
+}(W_{\tau_D}) \le\uppi\int h_{\nu^+}
(x) \sigma^+(\mathrm{d}x)= \bigl\langle\sigma^+ , \nu^+ \bigr\rangle_{\cM},
\]
and, after repeating the same argument, $\langle\sigma_D^-, \nu
^-\rangle_{\cM} \le\langle\sigma^-, \nu^-\rangle_{\cM}$ and so
on. Thus,
\[
\bigl\langle\sigma_D^+, \nu\bigr\rangle_{\cM} \le \bigl |\bigl
\langle\sigma^+ , \nu^+ \bigr\rangle_{\cM} \bigr | + \bigl |\bigl\langle\sigma^+ ,
\nu^- \bigr\rangle_{\cM} \bigr |,
\]
and a similar upper bound exists for $\langle\sigma_D^-, \nu\rangle
_{\cM}$. The above bounds imply
\[
\langle\sigma_D, \sigma_D \rangle_{\cM} \le \bigl |
\bigl\langle \sigma^+, \sigma^+ \bigr\rangle_{\cM}\bigr  | + 2\bigl  | \bigl\langle
\sigma ^+, \sigma^- \bigr\rangle_{\cM} \bigr | + \bigl | \bigl\langle\sigma^-,
\sigma^- \bigr\rangle_{\cM} \bigr |,
\]
which ensures that $\sigma_D \in\cM$. Now $D$ is a closed set with a
simply connected open interior and so we get
\[
|\sigma_D|\bigl(D^c\bigr) = \int\P_x(
W_{\tau_D} \notin D) |\sigma|(\dd x) =0.
\]
Therefore, $\sigma_D$ belongs to $\cM_D$.
Next, we establish that
\[
\E Z_\sigma Z_\nu= \E Z_{\sigma_D} Z_\nu\qquad
\forall \nu\in\cM_D.
\]
Since $\sigma_D$ is a measure that satisfies the relation
$\sigma_D(G) = \E_\sigma( 1_G(W_{\tau_D}))$,
for all Borel subsets $G$ of the plane, integrals with respect to
$\sigma_D$ can be defined as appropriate expected values of the
functions of $W_{\tau_D}$. In particular, for any element $f$ of an
appropriate class of functions, such an integral will satisfy
\[
\int f(x) \sigma_D(\dd x) = \E_\sigma f(W_{\tau_D}).
\]
Now take $f(x) = \E_x A_\nu(T)$. Then, the definition of the additive
function asserts that $\E Z_{\sigma_D} Z_\nu= \uppi\E_{\sigma_D}
A_\nu(T) $, but the above equation also implies
\[
\E_{\sigma_D} A_\nu(T) = \int\E_x
A_\nu(T) \sigma_D(\dd x) = \int f(x)
\sigma_D(\dd x) = \E_\sigma f(W_{\tau_D}) =
\E_\sigma \E_{W_{\tau_D}} A_\nu(T).
\]
Consequently, the strong Markov property of the Brownian motion in
equation (\ref{eq:strongmarkov}) applies, and we obtain
\[
E_{\sigma_D} A_\nu(T)=\E_\sigma A_\nu(
\tau_D+ T).
\]
Since $\tau_D$ is the first hitting time of $D$, the path of the
Brownian motion up to but not including time $\tau_D$ lies entirely
within the complement of $D$. However, the signed measure $\nu$
concentrates on $D$ making it imminent that $\E_\sigma A_\nu((0,\tau
_D)) =0$. And, therefore
\[
\cov( Z_{\sigma_D}, Z_\nu) = \uppi\E_{\sigma_D}
A_\nu(T) = \uppi\E _\sigma A_\nu(
\tau_D+ T)= \uppi\E_\sigma A_\nu(T) = \cov(
Z_{\sigma}, Z_\nu).
\]
This completes the proof.
\end{pf*}

Interestingly, when values are known along a straight line or on a
circle, the analytic formulas for the coefficient function $v(x,x_0)$
are available in closed form, making it possible to apply Theorem~\ref
{th:markov} directly to calculate relevant kriging predictions. For an
example, when $B$ is the unit circle, $v(x,x_0)$ becomes the Poisson kernel
%
\begin{equation}
\label{eq:poisson} v(x,x_0) = \frac{1}{2 \uppi} \frac{1- \|x_0\|^2}{\|x-x_0\|^2},\qquad  \|x\|
=1, \|x_0\|<1.
\end{equation}
Furthermore, when the boundary set $B$ is the $y$-axis, we refer to
\cite{r27} for a formula for the corresponding
coefficient
function $v(x,x_0)$.

Finally, we can also discuss predictions for functionals of the values
of the de Wijs process inside the boundary $B$ (e.g., $f(Z_\sigma)$
for some suitable function $f$), given $\{ Z_\nu\dvtx  \nu\in\cM_{D} \}
$, but this would require a knowledge of Wick products and Fock spaces
and is
beyond the scope of this paper.

\section{Discussion} \label{sec:discussion}

In practice, we only select finitely many regular or irregularly
distributed sampling locations and observe process values as aggregates
or averages over certain non-empty regular or irregular regions around
those sampling locations. In both instances there are certain
limitations in applying Matheron's kriging formula directly. For
example, if we observe only finitely many data values on the unit
circle, we won't be able to apply the exact kriging formula in (\ref
{eq:poisson}). Similarly, when the de Wijs process is used as a
statistical model for aggregates or averages of spatial variables over non-empty
regions in the two-dimensional plane, we simply lose the
Markov property because of aggregations or averaging. Examples include
agricultural field trials where
the variable of interest is the crop yield over plots, or disease
mapping where the spatial variable of interest is considered a
stochastically degraded version of an underlying unobserved spatial
component such as the log relative risk of non-infectious diseases over
a geographic region.

\begin{figure}
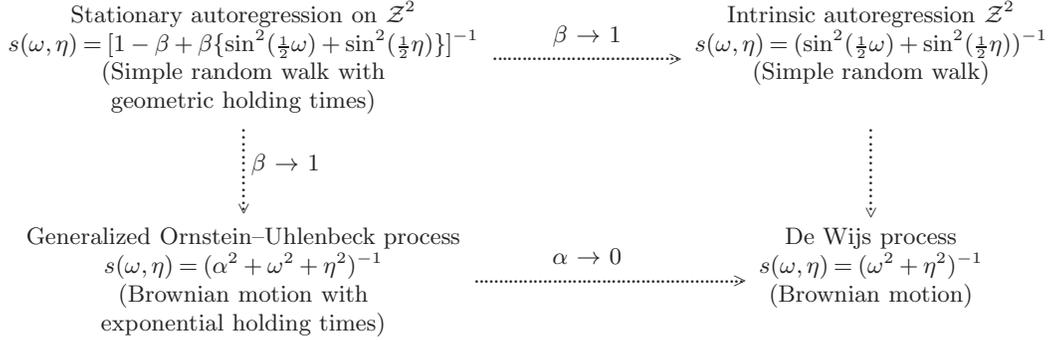

\begin{center}
{\fontsize{8.9pt}{10.9pt}\selectfont{\begin{diagram}
\begin{array}{c}
\mbox{Stationary autoregression on } \mathcal{Z}^2\\
s(\omega,\eta)=  [ 1-\beta+ \beta\{ \sin^2 (\shalf\omega) +
\sin^2( \shalf\eta)\} ]^{-1} \\
\mbox{(Simple random walk with}\\
\mbox{geometric holding times)}
\end{array}
&\rDotsto^{\hspace{.3in} \beta\, \rightarrow\, 1 \hspace{.3in}}
&
\begin{array}{c}
\mbox{Intrinsic autoregression } \mathcal{Z}^2\\
s(\omega,\eta) =  (\sin^2 (\shalf\omega) + \sin^2( \shalf
\eta) )^{-1} \\
\mbox{(Simple random walk)}\\
\mbox{}
\end{array}
\\
\\
\dDotsto_{\vspace{2in} \beta\, \rightarrow\, 1}
\\
&
&
\dDotsto_{}
\\
\begin{array}{c}
\mbox{Generalized Ornstein--Uhlenbeck process} \\
s(\omega, \eta)=(\alpha^2 + \omega^2 +\eta^2)^{-1} \\
\mbox{(Brownian motion with}\\
\mbox{exponential holding times)}
\end{array}
&
\rDotsto^{\alpha\, \rightarrow\, 0 }
&
\begin{array}{c}
\mbox{De Wijs process} \\
s(\omega,\eta)=(\omega^2 +\eta^2)^{-1} \\
\mbox{(Brownian motion)}\\
\mbox{}
\end{array}
\end{diagram}}}
\end{center}
\caption{Limit diagram for Gaussian Markov random fields and
associated Markov processes. Here $0 \le\beta< 1$ and $\alpha>0$.}
\label{fig:diagram}
\end{figure}

However, certain discrete approximations are possible, and, in fact, it
is the discrete approximations of the de Wijs process that have played
a major role in spatial statistics in the past thirty years; see, for
example, \cite{r2,r3,r22}, many subsequent papers,
and the books by Cressie \cite{r11}, Banerjee et al.
\cite{r1} and Rue and Held \cite{r45}. These discrete
approximations form a subclass of Gaussian
Markov random fields on regular and irregular lattices and have lattice
graph Laplacians as their precision (i.e., inverse of `covariance')
matrix. The nature of these approximations become clearer when we also
note that the log function (i.e., the generalized covariance of the de
Wijs process) is the inverse of the Laplacian on the plane.
Interestingly, these Gaussian random fields are also associated with
the random walks on the lattice graph, as the De Wijs process is with
the Browning motion on the plane. As a concrete example, the first
order symmetric intrinsic autoregressions on the two dimensional
integer lattice $\mathcal{Z}^2$ \cite{r5,r33}
is associated with the simple random walk on $\mathcal{Z}^2$. Thus
Dynkin's theory also applies here and we can obtain the coefficient
function of a corresponding kriging problem on the discrete lattice
$\mathcal{Z}^2$ from probabilities of the simple random walk at the
first hitting time. To summarize, the diagram in Figure~\ref
{fig:diagram} lists a few important spatial Gaussian Markov models and
the associated two-dimensional Markov processes. The top part of this
diagram notes the lattice Gaussian Markov random fields along with the
spectral densities and, in brackets, the associated lattice Markov
processes. The bottom part of the diagram provides limiting continuum
Gaussian Markov random fields with corresponding spectral densities
and, in brackets, the associated Markov processes. These continuum
random fields arise as the scaling limits of corresponding lattice
Markov fields from the top part of the diagram; see \cite{r36} and
\cite{r5} for details.

Another interesting point is that aggregates or averages of the de Wijs
process retain an approximate Markov
property that is known as the screening effect in geostatistics
\cite{r6,r49}. To give a very simple
example, we consider a regular lattice in the two-dimensional plane.
Let $X_{s,t}$ denote the average value of the spatial variable of
interest over the unit square whose center has integer coordinate
$(s,t)$ in the Euclidean plane. Given a realization $X^{S,T}$ of
$X_{s,t}$ for $s,t = -8, \ldots, 8$ but excluding $X_{0,0}$, we can
employ exact variogram computations \cite{r8,r33} to
find the coefficients $\omega_{s,t}$ in the conditional expectation
or ordinary kriging predictor \cite{r48}
%
\begin{equation}
\label{eq:kriging} \E\bigl( X_{0,0} | X^{S,T} \bigr) = \sum
\omega_{s,t} X_{s,t}
\end{equation}
under the regularized de Wijs process. The sum on the right-hand
extends over the aforementioned index set and the ordinary kriging
coefficients $\omega_{s,t}$ add up to 1. Table~\ref{tab:kriging}
shows the numerical values of $\omega_{s,t}$ to three decimals; in
view of symmetries only 44 of the $17^2 - 1 = 288$ coefficients need
to be shown. The screening effect is prominent here in that the
immediately neighboring cells dominate, with very few of the
remaining cells receiving non-negligible ordinary kriging
coefficients. However, it is not known to me if $\omega_{s,t}$ can be
interpreted in a meaningful way using certain probability calculations
of the Brownian motions, but here one can further try to derive
analytic form of $\omega_{s,t}$ for an infinite lattice from the
spectral density form of $X_{u,v}$. The same applies for aggregates or
averages of the first-order intrinsic autoregression. In general, it
would be interesting to know if one can better understand such an
approximate Markov property.

Some future directions can be added to this work. For example, the work
of \cite{r50} provides links between Dynkin's
isomorphisms and
constructions of statistical designs. Generalizations of Tjur's work in
the context of spatial designs would be an interesting matter for
future study.
%
\begin{table}
\tablewidth=\textwidth
\tabcolsep=0pt
\caption{Numerical values of the coefficients $\omega_{s,t}$ in the
ordinary kriging predictor (\protect\ref{eq:kriging}) under the
regularized de Wijs process} \label{tab:kriging}
\begin{tabular*}{\textwidth}{@{\extracolsep{\fill}}lllllllll@{}}
\hline
& \multicolumn{8}{l@{}}{$t$}\\[-5pt]
&\multicolumn{8}{l@{}}{\hrulefill}\\
$s$ & 1 & 2 & 3 & 4 & 5 & 6 & 7 & 8 \\
\hline
0 & \phantom{$-$}{0.342} & $-0.075$ & \phantom{$-$}{0.017} & $-0.004$ & 0.001 & 0.000 & 0.000 & 0.000
\\
1 & $-0.032$ & $-0.001$ & \phantom{$-$}{0.002} & $-0.001$ & 0.000 & 0.000 & 0.000 &
0.000 \\
2 & & \phantom{$-$}{0.002} & $-0.001$ & \phantom{$-$}{0.000} & 0.000 & 0.000 & 0.000 & 0.000 \\
3 & & & \phantom{$-$}{0.000} & \phantom{$-$}{0.000} & 0.000 & 0.000 & 0.000 & 0.000 \\
4 & & & & \phantom{$-$}{0.000} & 0.000 & 0.000 & 0.000 & 0.000 \\
5 & & & & & 0.000 & 0.000 & 0.000 & 0.000 \\
6 & & & & & & 0.000 & 0.000 & 0.000 \\
7 & & & & & & & 0.000 & 0.000 \\
8 & & & & & & & & 0.000 \\
\hline
\end{tabular*}
\end{table}
%

\section*{Acknowledgements}

Julian Besag drew my attention to the de Wijs process and suggested
that I study its Markov property. Tilmann Gneiting provided me with
edits and comments to an earlier version of the work. Atma Mandrekar
pointed me to the work of Dynkin during a visit to Michigan State
University. An anonymous reviewer, an associate editor and the editor
provided me with a careful reading of the paper. Finally, I acknowledge
support for the work from the National Science Foundation under award
DMS 0906300.


%

\printhistory


\begin{thebibliography}{57}

\bibitem{r1}
%
\begin{bbook}[auto:STB|2013/12/09|07:59:19]
\bauthor{\bsnm{Banerjee},~\bfnm{S.}\binits{S.}},
\bauthor{\bsnm{Carlin},~\bfnm{B.~P.}\binits{B.P.}} \AND
\bauthor{\bsnm{Gelfrand},~\bfnm{A.~E.}\binits{A.E.}}
(\byear{2004}).
\btitle{Hierarchical Modeling and Analysis for Spatial Data}.
\blocation{London}: \bpublisher{Chapman \& Hall}.
\bptok{imsref}%
\end{bbook}
%
\endbibitem

\bibitem{r2}
%
\begin{barticle}[mr]
\bauthor{\bsnm{Besag},~\bfnm{Julian}\binits{J.}}
(\byear{1974}).
\btitle{Spatial interaction and the statistical analysis of lattice systems}.
\bjournal{J. R. Stat. Soc. Ser.~B}
\bvolume{36}
\bpages{192--236}.
\bnote{With discussion by D. R. Cox, A. G. Hawkes, P. Clifford, P.
Whittle, K.
Ord, R. Mead, J. M. Hammersley and M. S. Bartlett and with a reply by the
author}.
\bid{issn={0035-9246}, mr={0373208}}
\bptnote{check related}%
\bptok{imsref}%
\end{barticle}
%
\endbibitem

\bibitem{r3}
%
\begin{barticle}[mr]
\bauthor{\bsnm{Besag},~\bfnm{Julian}\binits{J.}}
(\byear{1986}).
\btitle{On the statistical analysis of dirty pictures}.
\bjournal{J. R. Stat. Soc. Ser. B}
\bvolume{48}
\bpages{259--302}.
\bid{issn={0035-9246}, mr={0876840}}
\bptok{imsref}%
\end{barticle}
%
\endbibitem

\bibitem{r4}
%
\begin{barticle}[auto:STB|2013/12/09|07:59:19]
\bauthor{\bsnm{Besag},~\bfnm{J.}\binits{J.}}
(\byear{2002}).
\btitle{Discussion on the paper by McCullagh}.
\bjournal{Ann. Statist.}
\bvolume{30}
\bpages{1267--1277}.
\bptok{imsref}%
\end{barticle}
%
\endbibitem

\bibitem{r5}
%
\begin{barticle}[mr]
\bauthor{\bsnm{Besag},~\bfnm{Julian}\binits{J.}} \AND
\bauthor{\bsnm{Mondal},~\bfnm{Debashis}\binits{D.}}
(\byear{2005}).
\btitle{First-order intrinsic autoregressions and the de {W}ijs process}.
\bjournal{Biometrika}
\bvolume{92}
\bpages{909--920}.
\bid{doi={10.1093/biomet/92.4.909}, issn={0006-3444}, mr={2234194}}
\bptok{imsref}%
\end{barticle}
%
\endbibitem

\bibitem{r6}
%
\begin{bbook}[mr]
\bauthor{\bsnm{Chil{\`e}s},~\bfnm{Jean-Paul}\binits{J.P.}} \AND
\bauthor{\bsnm{Delfiner},~\bfnm{Pierre}\binits{P.}}
(\byear{1999}).
\btitle{Geostatistics: Modeling Spatial Uncertainty}.
\bseries{Wiley Series in Probability and Statistics: Applied
Probability and
Statistics}.
\blocation{New York}: \bpublisher{Wiley}.
\bid{doi={10.1002/9780470316993}, mr={1679557}}
\bptok{imsref}%
\end{bbook}
%
\endbibitem

\bibitem{r7}
%
\begin{bbook}[mr]
\bauthor{\bsnm{Chung},~\bfnm{Kai~Lai}\binits{K.L.}} \AND
\bauthor{\bsnm{Walsh},~\bfnm{John~B.}\binits{J.B.}}
(\byear{2005}).
\btitle{Markov Processes, {B}rownian Motion, and Time Symmetry},
\bedition{2nd}~ed.
\bseries{Grundlehren der Mathematischen Wissenschaften [Fundamental Principles
of Mathematical Sciences]}
\bvolume{249}.
\blocation{New York}: \bpublisher{Springer}.
\bid{mr={2152573}}
\bptok{imsref}%
\end{bbook}
%
\endbibitem

\bibitem{r8}
%
\begin{barticle}[mr]
\bauthor{\bsnm{Clifford},~\bfnm{David}\binits{D.}}
(\byear{2005}).
\btitle{Computation of spatial covariance matrices}.
\bjournal{J. Comput. Graph. Statist.}
\bvolume{14}
\bpages{155--167}.
\bid{doi={10.1198/106186005X27626}, issn={1061-8600}, mr={2137895}}
\bptok{imsref}%
\end{barticle}
%
\endbibitem

\bibitem{r9}
%
\begin{barticle}[mr]
\bauthor{\bsnm{Clifford},~\bfnm{David}\binits{D.}}
(\byear{2006}).
\btitle{Distribution of increases in residual log likelihood for
nested spatial
models}.
\bjournal{Comm. Statist. Simulation Comput.}
\bvolume{35}
\bpages{779--788}.
\bid{doi={10.1080/03610910600716886}, issn={0361-0918}, mr={2240044}}
\bptok{imsref}%
\end{barticle}
%
\endbibitem

\bibitem{r10}
%
\begin{barticle}[auto:STB|2013/12/09|07:59:19]
\bauthor{\bsnm{Clifford},~\bfnm{D.}\binits{D.}},
\bauthor{\bsnm{McBratney},~\bfnm{A.~B.}\binits{A.B.}},
\bauthor{\bsnm{Taylor},~\bfnm{J.}\binits{J.}} \AND
\bauthor{\bsnm{Whelan},~\bfnm{B.~M.}\binits{B.M.}}
(\byear{2006}).
\btitle{Generalized analysis of spatial variation in yield monitor data}.
\bjournal{J. Agric. Sci.}
\bvolume{144}
\bpages{45--51}.
\bptok{imsref}%
\end{barticle}
%
\endbibitem

\bibitem{r11}
%
\begin{bbook}[mr]
\bauthor{\bsnm{Cressie},~\bfnm{Noel A.~C.}\binits{N.A.C.}}
(\byear{1993}).
\btitle{Statistics for Spatial Data}.
\bseries{Wiley Series in Probability and Mathematical Statistics: Applied
Probability and Statistics}.
\blocation{New York}: \bpublisher{Wiley}.
\bnote{Revised reprint of the 1991 edition}.
\bid{mr={1239641}}
\bptok{imsref}%
\end{bbook}
%
\endbibitem

\bibitem{r12}
%
\begin{barticle}[mr]
\bauthor{\bsnm{Dobrushin},~\bfnm{R.~L.}\binits{R.L.}}
(\byear{1979}).
\btitle{Gaussian and their subordinated self-similar random generalized
fields}.
\bjournal{Ann. Probab.}
\bvolume{7}
\bpages{1--28}.
\bid{issn={0091-1798}, mr={0515810}}
\bptok{imsref}%
\end{barticle}
%
\endbibitem

\bibitem{n1}
%
\begin{bmisc}[auto:STB|2013/12/09|07:59:19]
\bauthor{\bsnm{Dutta},~\bfnm{S.}\binits{S.}} \AND
\bauthor{\bsnm{Mondal},~\bfnm{D.}\binits{D.}}
(\byear{2014}).
\bhowpublished{An $h$-likelihood method for spatial mixed linear models
based on intrinsic autoregressions. \textit{J. R. Stat. Soc. Ser. B Stat.
Method.} DOI:\doiurl{10.1111/rssb.12084}}.
\bptok{imsref}%
\end{bmisc}
%
\endbibitem

\bibitem{r13}
%
\begin{barticle}[mr]
\bauthor{\bsnm{Dynkin},~\bfnm{E.~B.}\binits{E.B.}}
(\byear{1980}).
\btitle{Markov processes and random fields}.
\bjournal{Bull. Amer. Math. Soc. (N.S.)}
\bvolume{3}
\bpages{975--999}.
\bid{doi={10.1090/S0273-0979-1980-14831-4}, issn={0273-0979}, mr={0585179}}
\bptok{imsref}%
\end{barticle}
%
\endbibitem

\bibitem{r14}
%
\begin{barticle}[mr]
\bauthor{\bsnm{Ekhaguere},~\bfnm{G.~O.~S.}\binits{G.O.S.}}
(\byear{1977}).
\btitle{On notions of {M}arkov property}.
\bjournal{J. Math. Phys.}
\bvolume{18}
\bpages{2104--2107}.
\bid{issn={0022-2488}, mr={0456149}}
\bptok{imsref}%
\end{barticle}
%
\endbibitem

\bibitem{r15}
%
\begin{bbook}[mr]
\bauthor{\bsnm{Gelfand},~\bfnm{I.~M.}\binits{I.M.}} \AND
\bauthor{\bsnm{Vilenkin},~\bfnm{N.~Ya.}\binits{N.Y.}}
(\byear{1964}).
\btitle{Generalized Functions. {V}ol. 4: Applications of Harmonic Analysis}.
\blocation{New York}: \bpublisher{Academic Press}.
\bnote{Translated from the Russian
by Amiel
Feinstein}.
\bid{mr={0435834}}
\bptnote{check year}%
\bptok{imsref}%
\end{bbook}
%
\endbibitem

\bibitem{r16}
%
\begin{bbook}[mr]
\bauthor{\bsnm{Glimm},~\bfnm{James}\binits{J.}} \AND
\bauthor{\bsnm{Jaffe},~\bfnm{Arthur}\binits{A.}}
(\byear{1981}).
\btitle{Quantum Physics: A Functional Integral Point of View}.
\blocation{New York}: \bpublisher{Springer}.
\bid{mr={0628000}}
\bptok{imsref}%
\end{bbook}
%
\endbibitem

\bibitem{r17}
%
\begin{barticle}[mr]
\bauthor{\bsnm{Holley},~\bfnm{R.}\binits{R.}} \AND
\bauthor{\bsnm{Stroock},~\bfnm{D.}\binits{D.}}
(\byear{1980}).
\btitle{The {D}.{L}.{R}. conditions for translation invariant {G}aussian
measures on $\mathcal{S}^{\prime} (\mathcal{R}^{d})$}.
\bjournal{Z. Wahrsch. Verw. Gebiete}
\bvolume{53}
\bpages{293--304}.
\bid{doi={10.1007/BF00531439}, issn={0044-3719}, mr={0586022}}
\bptok{imsref}%
\end{barticle}
%
\endbibitem

\bibitem{r18}
%
\begin{bbook}[auto:STB|2013/12/09|07:59:19]
\bauthor{\bsnm{Journel},~\bfnm{A.~G.}\binits{A.G.}} \AND
\bauthor{\bsnm{Huijbregts},~\bfnm{C.~H.}\binits{C.H.}}
(\byear{1978}).
\btitle{Mining Geostatistics}.
\blocation{London}: \bpublisher{Academic Press}.
\bptok{imsref}%
\end{bbook}
%
\endbibitem

\bibitem{r19}
%
\begin{barticle}[mr]
\bauthor{\bsnm{Kallianpur},~\bfnm{G.}\binits{G.}} \AND
\bauthor{\bsnm{Mandrekar},~\bfnm{V.}\binits{V.}}
(\byear{1974}).
\btitle{The {M}arkov property for generalized {G}aussian random fields}.
\bjournal{Ann. Inst. Fourier (Grenoble)}
\bvolume{24}
\bpages{143--167}.
\bid{issn={0373-0956}, mr={0405569}}
\bptok{imsref}%
\end{barticle}
%
\endbibitem

\bibitem{r20}
%
\begin{barticle}[mr]
\bauthor{\bsnm{Kenyon},~\bfnm{Richard}\binits{R.}}
(\byear{2001}).
\btitle{Dominos and the {G}aussian free field}.
\bjournal{Ann. Probab.}
\bvolume{29}
\bpages{1128--1137}.
\bid{doi={10.1214/aop/1015345599}, issn={0091-1798}, mr={1872739}}
\bptok{imsref}%
\end{barticle}
%
\endbibitem

\bibitem{r21}
%
\begin{barticle}[mr]
\bauthor{\bsnm{K{\"u}nsch},~\bfnm{H.}\binits{H.}}
(\byear{1979}).
\btitle{Gaussian {M}arkov random fields}.
\bjournal{J. Fac. Sci. Univ. Tokyo Sect. IA Math.}
\bvolume{26}
\bpages{53--73}.
\bid{issn={0040-8980}, mr={0539773}}
\bptok{imsref}%
\end{barticle}
%
\endbibitem

\bibitem{r22}
%
\begin{barticle}[mr]
\bauthor{\bsnm{K{\"u}nsch},~\bfnm{Hans~R.}\binits{H.R.}}
(\byear{1987}).
\btitle{Intrinsic autoregressions and related models on the two-dimensional
lattice}.
\bjournal{Biometrika}
\bvolume{74}
\bpages{517--524}.
\bid{issn={0006-3444}, mr={0909356}}
\bptok{imsref}%
\end{barticle}
%
\endbibitem


\bibitem{r24}
%
\begin{binproceedings}[mr]
\bauthor{\bsnm{L{\'e}vy},~\bfnm{Paul}\binits{P.}}
(\byear{1956}).
\btitle{A special problem of {B}rownian motion, and a general theory of
{G}aussian random functions}.
In \bbooktitle{Proceedings of the {T}hird {B}erkeley {S}ymposium on
{M}athematical {S}tatistics and {P}robability, 1954--1955, Vol. {II}}
\bpages{133--175}.
\blocation{Berkeley}: \bpublisher{Univ. California Press}.
\bid{mr={0090934}}
\bptok{imsref}%
\end{binproceedings}
%
\endbibitem

\bibitem{r25}
%
\begin{bincollection}[mr]
\bauthor{\bsnm{Mandrekar},~\bfnm{V.}\binits{V.}} \AND
\bauthor{\bsnm{Zhang},~\bfnm{Sixiang}\binits{S.}}
(\byear{1993}).
\btitle{Markov property of measure-indexed {G}aussian random fields}.
In \bbooktitle{Stochastic Processes}
\bpages{253--262}.
\blocation{New York}: \bpublisher{Springer}.
\bid{mr={1427321}}
\bptok{imsref}%
\end{bincollection}
%
\endbibitem


\bibitem{r26}
%
\begin{bbook}[auto:STB|2013/12/09|07:59:19]
\bauthor{\bsnm{Matheron},~\bfnm{G.}\binits{G.}}
(\byear{1962}).
\btitle{Trait\'e de G\'eostatistique Appliqu\'ee}
\bvolume{1}.
\blocation{Paris}: \bpublisher{Editions Technip}.
\bptok{imsref}%
\end{bbook}
%
\endbibitem

\bibitem{r27}
%
\begin{bmisc}[auto:STB|2013/12/09|07:59:19]
\bauthor{\bsnm{Matheron},~\bfnm{G.}\binits{G.}}
(\byear{1963}).
\bhowpublished{\textit{Processus Markoviens Normaux Stationaires \`a $n$ Dimensions. Note
G\'eostatistique} \textbf{50}. Paris: Centre de G\'eostatistique, Ecole des Mines de
 Paris}.
\bptok{imsref}%
\end{bmisc}
%
\endbibitem

\bibitem{r27a}
%
\begin{bmisc}[auto:STB|2013/12/09|07:59:19]
\bauthor{\bsnm{Matheron},~\bfnm{G.}\binits{G.}}
(\byear{1971}).
\bhowpublished{\textit{The Theory of Regionalized Variables and Its
Applications. Les Cahiers du Centre de Morphologie Math\'ematique de
Fontainebleau} \textbf{5}}.
\bptok{imsref}%
\end{bmisc}
%
\endbibitem


\bibitem{r28}
%
\begin{barticle}[mr]
\bauthor{\bsnm{Mattner},~\bfnm{L.}\binits{L.}}
(\byear{1997}).
\btitle{Strict definiteness of integrals via complete monotonicity of
derivatives}.
\bjournal{Trans. Amer. Math. Soc.}
\bvolume{349}
\bpages{3321--3342}.
\bid{doi={10.1090/S0002-9947-97-01966-1}, issn={0002-9947}, mr={1422615}}
\bptok{imsref}%
\end{barticle}
%
\endbibitem

\bibitem{r29}
%
\begin{barticle}[mr]
\bauthor{\bsnm{McCullagh},~\bfnm{Peter}\binits{P.}}
(\byear{2002}).
\btitle{What is a statistical model?}
\bjournal{Ann. Statist.}
\bvolume{30}
\bpages{1225--1310}.
\bnote{With comments and a rejoinder by the author}.
\bid{doi={10.1214/aos/1035844977}, issn={0090-5364}, mr={1936320}}
\bptnote{check related}%
\bptok{imsref}%
\end{barticle}
%
\endbibitem

\bibitem{r30}
%
\begin{barticle}[auto:STB|2013/12/09|07:59:19]
\bauthor{\bsnm{McCullagh},~\bfnm{P.}\binits{P.}} \AND
\bauthor{\bsnm{Clifford},~\bfnm{D.}\binits{D.}}
(\byear{2006}).
\btitle{Evidence for conformal invariance of crop yields}.
\bjournal{Proc. R. Soc. Lond. Ser. A Math. Phys. Eng. Sci.}
\bvolume{462}
\bpages{2119--2143}.
\bptok{imsref}%
\end{barticle}
%
\endbibitem

\bibitem{r31}
%
\begin{barticle}[mr]
\bauthor{\bsnm{McKean},~\bfnm{H.~P.}\binits{H.P.} \bsuffix{Jr.}}
(\byear{1963}).
\btitle{Brownian motion with a several-dimensional time}.
\bjournal{Theory Probab. Appl.}
\bvolume{8}
\bpages{335--354}.
\bid{issn={0040-361X}}
\bptok{imsref}%
\end{barticle}
%
\endbibitem

\bibitem{r32}
%
\begin{barticle}[auto:STB|2013/12/09|07:59:19]
\bauthor{\bsnm{Molchan},~\bfnm{G.~M.}\binits{G.M.}}
(\byear{1971}).
\btitle{Characterization of Gaussian fields with Markov property}.
\bjournal{Soviet Math. Dokl.}
\bvolume{12}
\bpages{563--567}.
\bptok{imsref}%
\end{barticle}
%
\endbibitem

\bibitem{r33}
%
\begin{bmisc}[auto:STB|2013/12/09|07:59:19]
\bauthor{\bsnm{Mondal},~\bfnm{D.}\binits{D.}}
(\byear{2005}).
\bhowpublished{Variogram calculations for first-order intrinsic
autoregressions and
the de Wijs process. Technical Report 479, Department of Statistics,
 University of Washington}.
\bptok{imsref}%
\end{bmisc}
%
\endbibitem

\bibitem{r34}
%
\begin{bmisc}[auto:STB|2013/12/09|07:59:19]
\bauthor{\bsnm{Mondal},~\bfnm{D.}\binits{D.}}
(\byear{2015}).
\bhowpublished{Generalized Gaussian Markov random fields and modeling
disease risk. Under revision}.
\bptok{imsref}%
\end{bmisc}
%
\endbibitem



\bibitem{r35}
%
\begin{bmisc}[auto:STB|2013/12/09|07:59:19]
\bauthor{\bsnm{Mondal},~\bfnm{D.}\binits{D.}}
(\byear{2015}).
\bhowpublished{On Tobler's pycnophylactic interpolation. Unpublished manuscript,
Oregon State Univ.}
\bptok{imsref}%
\end{bmisc}
%
\endbibitem

\bibitem{r36}
%
\begin{barticle}[mr]
\bauthor{\bsnm{Moran},~\bfnm{P.~A.~P.}\binits{P.A.P.}}
(\byear{1973}).
\btitle{A {G}aussian {M}arkovian process on a square lattice}.
\bjournal{J. Appl. Probab.}
\bvolume{10}
\bpages{54--62}.
\bid{issn={0021-9002}, mr={0353437}}
\bptok{imsref}%
\end{barticle}
%
\endbibitem

\bibitem{r37}
%
\begin{barticle}[mr]
\bauthor{\bsnm{Nelson},~\bfnm{Edward}\binits{E.}}
(\byear{1973}a).
\btitle{Construction of quantum fields from {M}arkoff fields}.
\bjournal{J. Funct. Anal.}
\bvolume{12}
\bpages{97--112}.
\bid{mr={0343815}}
\bptok{imsref}%
\end{barticle}
%
\endbibitem

\bibitem{r38}
%
\begin{barticle}[mr]
\bauthor{\bsnm{Nelson},~\bfnm{Edward}\binits{E.}}
(\byear{1973}b).
\btitle{The free {M}arkoff field}.
\bjournal{J. Funct. Anal.}
\bvolume{12}
\bpages{211--227}.
\bid{mr={0343816}}
\bptok{imsref}%
\end{barticle}
%
\endbibitem

\bibitem{r39}
%
\begin{bbook}[mr]
\bauthor{\bsnm{Port},~\bfnm{Sidney~C.}\binits{S.C.}} \AND
\bauthor{\bsnm{Stone},~\bfnm{Charles~J.}\binits{C.J.}}
(\byear{1978}).
\btitle{Brownian Motion and Classical Potential Theory}.
\bseries{Probability and Mathematical Statistics}.
\blocation{New York}: \bpublisher{Academic Press}.
\bid{mr={0492329}}
\bptok{imsref}%
\end{bbook}
%
\endbibitem

\bibitem{r40}
%
\begin{bbook}[mr]
\bauthor{\bsnm{Rao},~\bfnm{Murali}\binits{M.}}
(\byear{1977}).
\btitle{Brownian Motion and Classical Potential Theory}.
\bseries{Lecture Notes Series}
\bvolume{47}.
\blocation{Aarhus}: \bpublisher{Matematisk Institut, Aarhus Univ.}
\bid{mr={0440718}}
\bptok{imsref}%
\end{bbook}
%
\endbibitem

\bibitem{r41}
%
\begin{barticle}[mr]
\bauthor{\bsnm{R{\"o}ckner},~\bfnm{Michael}\binits{M.}}
(\byear{1983}).
\btitle{Markov property of generalized fields and axiomatic potential theory}.
\bjournal{Math. Ann.}
\bvolume{264}
\bpages{153--177}.
\bid{doi={10.1007/BF01457522}, issn={0025-5831}, mr={0711875}}
\bptok{imsref}%
\end{barticle}
%
\endbibitem

\bibitem{r42}
%
\begin{barticle}[mr]
\bauthor{\bsnm{R{\"o}ckner},~\bfnm{Michael}\binits{M.}}
(\byear{1985}).
\btitle{Generalized {M}arkov fields and {D}irichlet forms}.
\bjournal{Acta Appl. Math.}
\bvolume{3}
\bpages{285--311}.
\bid{doi={10.1007/BF00047332}, issn={0167-8019}, mr={0790552}}
\bptok{imsref}%
\end{barticle}
%
\endbibitem

\bibitem{r43}
%
\begin{barticle}[auto:STB|2013/12/09|07:59:19]
\bauthor{\bsnm{Rozanov},~\bfnm{Y.~A.}\binits{Y.A.}}
(\byear{1977}).
\btitle{Markovian random fields and stochastic partial differential equations}.
\bjournal{Mat. USSR--Sb.}
\bvolume{32}
\bpages{515--534}.
\bptok{imsref}%
\end{barticle}
%
\endbibitem

\bibitem{r44}
%
\begin{bincollection}[mr]
\bauthor{\bsnm{Rozanov},~\bfnm{Yu.~A.}\binits{Y.A.}}
(\byear{1979}).
\btitle{Stochastic {M}arkovian fields}.
In \bbooktitle{Developments in Statistics, {V}ol. 2}
\bpages{203--234}.
\blocation{New York}: \bpublisher{Academic Press}.
\bid{mr={0554181}}
\bptok{imsref}%
\end{bincollection}
%
\endbibitem

\bibitem{r45}
%
\begin{bbook}[mr]
\bauthor{\bsnm{Rue},~\bfnm{H{\aa}vard}\binits{H.}} \AND
\bauthor{\bsnm{Held},~\bfnm{Leonhard}\binits{L.}}
(\byear{2005}).
\btitle{Gaussian {M}arkov Random Fields: Theory and Applications}.
\bseries{Monographs on Statistics and Applied Probability}
\bvolume{104}.
\blocation{Boca Raton, FL}: \bpublisher{Chapman \& Hall/CRC}.
\bid{doi={10.1201/9780203492024}, mr={2130347}}
\bptok{imsref}%
\end{bbook}
%
\endbibitem

\bibitem{r46}
%
\begin{barticle}[mr]
\bauthor{\bsnm{Sch{\"a}fer},~\bfnm{J{\"o}rg}\binits{J.}}
(\byear{1996}).
\btitle{Abstract {M}arkov property and local operators}.
\bjournal{J. Funct. Anal.}
\bvolume{138}
\bpages{137--169}.
\bid{doi={10.1006/jfan.1996.0059}, issn={0022-1236}, mr={1391633}}
\bptok{imsref}%
\end{barticle}
%
\endbibitem

\bibitem{r47}
%
\begin{barticle}[mr]
\bauthor{\bsnm{Sheffield},~\bfnm{Scott}\binits{S.}}
(\byear{2007}).
\btitle{Gaussian free fields for mathematicians}.
\bjournal{Probab. Theory Related Fields}
\bvolume{139}
\bpages{521--541}.
\bid{doi={10.1007/s00440-006-0050-1}, issn={0178-8051}, mr={2322706}}
\bptok{imsref}%
\end{barticle}
%
\endbibitem

\bibitem{r48}
%
\begin{bbook}[mr]
\bauthor{\bsnm{Stein},~\bfnm{Michael~L.}\binits{M.L.}}
(\byear{1999}).
\btitle{Interpolation of Spatial Data: Some Theory for Kriging}.
\bseries{Springer Series in Statistics}.
\blocation{New York}: \bpublisher{Springer}.
\bid{doi={10.1007/978-1-4612-1494-6}, mr={1697409}}
\bptok{imsref}%
\end{bbook}
%
\endbibitem

\bibitem{r49}
%
\begin{barticle}[mr]
\bauthor{\bsnm{Stein},~\bfnm{Michael~L.}\binits{M.L.}}
(\byear{2002}).
\btitle{The screening effect in kriging}.
\bjournal{Ann. Statist.}
\bvolume{30}
\bpages{298--323}.
\bid{doi={10.1214/aos/1015362194}, issn={0090-5364}, mr={1892665}}
\bptok{imsref}%
\end{barticle}
%
\endbibitem

\bibitem{r50}
%
\begin{barticle}[mr]
\bauthor{\bsnm{Tjur},~\bfnm{Tue}\binits{T.}}
(\byear{1991}).
\btitle{Block designs and electrical networks}.
\bjournal{Ann. Statist.}
\bvolume{19}
\bpages{1010--1027}.
\bid{doi={10.1214/aos/1176348134}, issn={0090-5364}, mr={1105858}}
\bptok{imsref}%
\end{barticle}
%
\endbibitem

\bibitem{r51}
%
\begin{barticle}[mr]
\bauthor{\bsnm{Tobler},~\bfnm{Waldo~R.}\binits{W.R.}}
(\byear{1979}).
\btitle{Smooth pycnophylactic interpolation for geographical regions}.
\bjournal{J. Amer. Statist. Assoc.}
\bvolume{74}
\bpages{519--536}.
\bnote{With a comment by Nira Dyn [Nira Richter-Dyn], Grace Wahba and
Wing Hung
Wong and a rejoinder by the author}.
\bid{issn={0003-1291}, mr={0548256}}
\bptnote{check related}%
\bptok{imsref}%
\end{barticle}
%
\endbibitem

\bibitem{r52}
%
\begin{barticle}[mr]
\bauthor{\bsnm{Urbanik},~\bfnm{K.}\binits{K.}}
(\byear{1962}).
\btitle{Generalized stationary processes of {M}arkovian character}.
\bjournal{Studia Math.}
\bvolume{21}
\bpages{261--282}.
\bid{issn={0039-3223}, mr={0150835}}
\bptnote{check year}%
\bptok{imsref}%
\end{barticle}
%
\endbibitem

\bibitem{r53}
%
\begin{barticle}[mr]
\bauthor{\bsnm{Wong},~\bfnm{E.}\binits{E.}}
(\byear{1969}).
\btitle{Homogeneous {G}auss--{M}arkov random fields}.
\bjournal{Ann. Math. Statist.}
\bvolume{40}
\bpages{1625--1634}.
\bid{issn={0003-4851}, mr={0263148}}
\bptok{imsref}%
\end{barticle}
%
\endbibitem

\bibitem{r54}
%
\begin{barticle}[mr]
\bauthor{\bsnm{Yaglom},~\bfnm{A.~M.}\binits{A.M.}}
(\byear{1957}).
\btitle{Certain types of random fields in {$n$}-dimensional space
similar to
stationary stochastic processes}.
\bjournal{Theory Probab. Appl.}
\bvolume{2}
\bpages{273--320}.
\bid{issn={0040-361X}}
\bptok{imsref}%
\end{barticle}
%
\endbibitem

\end{thebibliography}
\end{document}